\documentclass[12pt,leqno]{article}

\newcommand{\EulerSol}{%
\unitlength 0.1in%
\ifmmode%
\begin{picture}(  1.2000,  1.2000)(  0.0000, -1.1000)
\special{pn 8}%
\special{ar 60 60 60 60  0.0000000 6.2831853}%
\special{pn 20}%
\special{sh 1}%
\special{ar 60 60 10 10 0  6.28318530717959E+0000}%
\special{sh 1}%
\special{ar 60 60 10 10 0  6.28318530717959E+0000}%
\end{picture}%
\else%
\begin{picture}(  1.8000,  1.2000)(  0.0000, -1.0000)
\special{pn 8}%
\special{ar 60 60 60 60  0.0000000 6.2831853}%
\special{pn 20}%
\special{sh 1}%
\special{ar 60 60 10 10 0  6.28318530717959E+0000}%
\special{sh 1}%
\special{ar 60 60 10 10 0  6.28318530717959E+0000}%
\end{picture}%
\fi}

\newcommand{\EulerLuna}{%
\unitlength 0.1in%
\ifmmode%
\begin{picture}(  0.8000,  1.2000)(  0.0000, -1.2000)
\special{pn 8}%
\special{ar 20 60 60 60  4.7123890 6.2831853}%
\special{ar 20 60 60 60  0.0000000 1.5707963}%
\special{pn 8}%
\special{pa 20 120}%
\special{pa 0 120}%
\special{fp}%
\special{pn 8}%
\special{pa 20 0}%
\special{pa 0 0}%
\special{fp}%
\special{pn 8}%
\special{ar -20 60 60 60  5.1760366 6.2831853}%
\special{ar -20 60 60 60  0.0000000 1.2036225}%
\end{picture}%
\else%
\begin{picture}(  1.6000,  1.2000)(  -0.1000, -1.0000)
\special{pn 8}%
\special{ar 20 60 60 60  4.7123890 6.2831853}%
\special{ar 20 60 60 60  0.0000000 1.5707963}%
\special{pn 8}%
\special{pa 20 120}%
\special{pa 0 120}%
\special{fp}%
\special{pn 8}%
\special{pa 20 0}%
\special{pa 0 0}%
\special{fp}%
\special{pn 8}%
\special{ar -20 60 60 60  5.1760366 6.2831853}%
\special{ar -20 60 60 60  0.0000000 1.2036225}%
\end{picture}
\fi}

\usepackage{amsfonts}
\usepackage{amsmath}
\usepackage{amssymb}

\newcommand{\de}{\delta}
\newcommand{\G}{\Gamma}
\newcommand{\qed }{\hfill\hbox{\rule{5pt}{8pt}}\medskip}
\newcommand{\ra}{\rightarrow}
\newcommand{\re}{\mbox{Re}}
\newcommand{\wt}{\widetilde}
\newcommand{\z}{\zeta}

\newcommand{\C}{{\bf C}}
\newcommand{\Z}{{\bf Z}}

\begin{document}

\begin{center}
  {\large \bf A variation of Euler's approach to values of \\
    the Riemann zeta function}
  \vspace{2.5ex}\\
  {\sc Masanobu Kaneko, Nobushige Kurokawa \\
    and Masato Wakayama\footnotetext{Work in part supported by
      Grant-in Aid for Scientific Research (B) No.11440010, and by
      Grant-in Aid for Exploratory Research No.13874004, Japan
      Society for the Promotion of Science.}}\\
\end{center}
\begin{center}
{\it Dedicated to Leonhard Euler on his 296th birthday}
\end{center}



\begin{center}
{\bf 1. Introduction}
\end{center}

Sometime around 1740, Euler \cite{E3,E4,E5} took a decisive step to
unraveling the  functional equation of the Riemann zeta function when
he discovered  a marvelous method of calculating the values of the
(absolutely!) divergent series 
\begin{eqnarray*}
``1+1+1+1+1+\cdots" &=& -\frac12,\\
``1+2+3+4+5+\cdots" &=& -\frac1{12},\\
``1+4+9+16+25+\cdots" &=&0, \\
``1+8+27+64+125+\cdots" &=&\frac1{120},\; \text{etc.}.
\end{eqnarray*}
In modern terms, these are the values at non-positive integer
arguments of the Riemann zeta function $\zeta(s)$, which is defined by the
series, absolutely convergent in $\re(s)>1$:
$$
\zeta(s)
=1+\frac1{2^s}+\frac1{3^s}+\frac1{4^s}+\frac1{5^s}+\cdots.
$$
Naturally, having no notion of the analytic continuation at that
time,  to say nothing of functions of complex variable, Euler had to
find a way of giving a meaning to those values of divergent series.
What he actually did proceeds as follows.  First, he directs his
attention to ``less divergent'' alternating  series
$$
\EulerSol \quad \quad \quad
1^m-2^m+3^m-4^m+5^m-6^m+7^m-8^m+\text{etc.}
$$
since its convergent counterpart
$$
\EulerLuna \quad \quad \quad
\dfrac1{1^n}-\dfrac1{2^n}+\dfrac1{3^n}-\dfrac1{4^n}%
+\dfrac1{5^n}-\dfrac1{6^n}+\dfrac1{7^n}-\dfrac1{8^n}+\text{etc.}
$$
does indeed have faster convergence and is linked with the
original series  by the simple relation
\begin{equation}
  \label{eq:alt}
\wt\zeta(s)
=(1-2^{1-s})\zeta(s),
\end{equation}
where we have put 
$$
\wt\zeta(s)
=1-\frac1{2^s}+\frac1{3^s}-\frac1{4^s}+\frac1{5^s}-\cdots.
$$
Then, he observes that the value \EulerSol
is obtained as a ``limit'' of the power series 
\begin{equation}
  \label{eq:eulpow}
  1^m-2^mx+3^mx^2-4^mx^3+5^mx^4-6^mx^5+\text{etc.}
\end{equation}
as $x\ra1$, since, although the series itself converges only for
$|x|<1$, it has an expression as a rational function (analytic
continuation, as we now put it), finite at $x=1$, which is obtained by
a successive application of multiplication by $x$ and differentiation
(or equivalently, applying the {\it Euler operator} $x\cdot d/dx$
successively after once multiplied by $x$) to the geometric series
expansion
\begin{equation}
  \label{eq:geom}
  \frac1{1+x}=1-x+x^2-x^3+x^4-x^5+\cdots, \quad (|x|<1). 
\end{equation}
For instance, if we substitute $x=1$ in \eqref{eq:geom}, we find
formally
$$
\frac12=1-1+1-1-1-\cdots=\wt\zeta(0)
$$
and hence, in view of \eqref{eq:alt},  we have $\zeta(0)=-1/2$. 
A few more examples are
\begin{eqnarray*}
  \frac1{(1+x)^2}&=&1-2x+3x^2-4x^3+5x^4-\cdots, \\
\frac{1-x}{(1+x)^3}&=&1-2^2x+3^2x^2-4^2x^3+5^2x^4-\cdots,\\
\frac{1-4x+x^2}{(1+x)^4}&=&1-2^3x+3^3x^2-4^3x^3+5^3x^4-\cdots,\;
\text{etc.}.
\end{eqnarray*}
which give us
$$
\wt\zeta(-1)=\frac14,\quad \wt\zeta(-2)=0,\quad
\wt\zeta(-3)=-\frac18,\cdots
$$
and in turn
$$
 \zeta(-1)=-\frac1{12},\quad  \zeta(-2)=0,\quad
 \zeta(-3)=\frac1{120},\cdots.
$$

For all that splendid idea however, this method unfortunately provides
no rigorous way to establish the values of $\z(-m)$ as values of the
analytically continued function $\z(s)$ at $s=-m$.   (This may only be
an afterthought and merely shows the degree to which Euler was ahead of
his time and how much our modern point of view owes to him.)

In the present article, aiming to evaluate the value $\z(-m)$ in as
elementary, and yet rigorous, way as possible as the value of
analytically continued function $\z(s)$, we introduce and investigate
a new $q$-analogue of the Riemann zeta function.  As becomes clear in
the course of our study, this function serves very well for the
purpose not only of computing $\z(-m)$ but  also of giving a nice
$q$-analogue of $\z(s)$ valid for {\em all} $s\in\C$.

To be more specific, as an alternative for a series like
\eqref{eq:eulpow}, we put $x=q^t$ (note $x\cdot d/dx$ is essentially
$d/dt$) and,  instead of repeating differentiation (this inevitably
restricts us to looking only at the integer arguments), we replace
$n^m$ by the $q$-integer $[n]_q:=(1-q^n)/(1-q)$ raised by the power
$-s$ (recall that Euler is the grand ``Master of $q$''!); namely, we
consider the series
\begin{equation}
  \label{eq:fqst}
  f_q(s,t):=\sum_{n=1}^\infty \frac{q^{nt}}{[n]_q^s}
=\frac{q^t}{[1]_q^s}+\frac{q^{2t}}{[2]_q^s}+\frac{q^{3t}}{[3]_q^s}
+\frac{q^{4t}}{[4]_q^s}+\cdots.
\end{equation}
Throughout the paper, we always assume $0<q<1$, so the series
\eqref{eq:fqst} converges absolutely for any $s\in\C$ and $\re(t)>0$.
If $\re(s)>1$ (and $\re(t)>0$), the series obviously converges to
$\z(s)$ when $q\uparrow1$. This suggests that we should regard the
function $f_q(s,t)$ as a $q$-analogue of the Riemann zeta function
$\z(s)$, but we reserve this until we make the specialization $t=s-1$
which turns out to be utterly crucial. Before going into the
specialization, we establish  below the meromorphic continuation of
$f_q(s,t)$ as a function of two variables $s$ and $t$, which is
carried out quite easily by using the binomial theorem.

In the next section, we specialize $t=s-1$  and establish a formula
for $s=-m\in\Z_{\le0}$ (Proposition 2) as well as its limit when
$q\uparrow1$ (Theorem 1). Then we give the result concerning  the
limit as $q\uparrow1$ for any $s$ (Theorem 2).

\noindent {\bf Proposition 1.} \begin{em}Let $0<q<1$. 
  As a function of $(s,t)\in\C^2$, $f_q(s,t)$ is continued
  meromorphically via the series expansion
\begin{eqnarray*}
f_q(s,t)&=&(1-q)^{s}\sum_{r=0}^\infty \binom{s+r-1}{r}
\frac{q^{t+r}}{1-q^{t+r}}\\
&=& (1-q)^s\left(\frac{q^t}{1-q^{t}}+s\frac{q^{t+1}}{1-q^{t+1}}
+\frac{s(s+1)}2\frac{q^{t+2}}{1-q^{t+2}}
+\cdots\right),
\end{eqnarray*}
having poles of order $1$ at $t\in\Z_{\le0}+2\pi i\Z/\log q
:=\left\{a +2\pi i b/\log q\,\vert\, a,b\in\Z,\,
  a\le0\right\}$.
\end{em}

\noindent {\it Proof.} ~We just apply the binomial 
expansion $(1-q^{n})^{-s}=\sum_{r= 0}^\infty\binom{s+r-1}{r}q^{nr}$
and change the order of summations to get
\begin{eqnarray*}
f_q(s,t)&=&(1-q)^s\sum_{n=1}^\infty \frac{q^{nt}}{(1-q^{n})^{s}}
=(1-q)^s\sum_{n=1}^\infty q^{nt}
\sum_{r=0}^\infty \binom{s+r-1}{r}q^{nr}\\
&=&(1-q)^s\sum_{r=0}^\infty \binom{s+r-1}{r}\sum_{n=1}^\infty
q^{n(t+r)}=(1-q)^s\sum_{r=0}^\infty \binom{s+r-1}{r}\frac{q^{t+r}}{1-q^{t+r}}.
\end{eqnarray*}
The other assertions follows readily from this. \qed

\noindent
{\it Remark.} It is worth noting that the function $f_q(s,t)$ can be
expressed as the (beta-like) Jackson integral. In fact, we have
$$
q^{-t}(1-q)^{1-s}f_q(s,t) =(1-q)\sum_{j=0}^\infty
\frac{q^{jt}}{(1-q^{j+1})^s}  =\int_0^1 x^{t-1}(1-qx)^{-s}d_qx.
$$

\begin{center}
{\bf  2. Main results}
\end{center}

Now we put $t=s-1$. When $s=-m\in\Z_{\le0}$, the point
$(s,t)=(-m,-m-1)$ lies on the pole divisor $t=-m-1$ of
$f_q(s,t)$. Nevertheless, a sort of ``miracle''  happens that the
point turns out to be  what is called ``the point of indeterminacy'',
the function $f_q(s,s-1)$ having a finite limit as $s\ra -m$  and
moreover the limit approaches to the ``correct'' value $\zeta(-m)$ as
$q\uparrow1$. What is more, the function $f_q(s,s-1)$ converges as
$q\uparrow1$ to $\z(s)$ for {\it any} $s\,$!  These results, to be
proved in quite elementary ways (certainly with only devices of which
Euler could avail himself), well justify the function $f_q(s,s-1)$
being referred to as the ``true'' $q$-analogue of the Riemann zeta
function, and we label it hereafter as
$$
\z_q(s):=f_q(s,s-1)=\sum_{n=1}^\infty \frac{q^{n(s-1)}}{[n]_q^s}
=\frac{q^{s-1}}{[1]_q^s}+\frac{q^{2(s-1)}}{[2]_q^s}
+\frac{q^{3(s-1)}}{[3]_q^s} +\frac{q^{4(s-1)}}{[4]_q^s}+\cdots.$$

\noindent
{\it Remark.} 1) Proper choice of $t$ seems to be essential. For
example, the choice $t=s$ adopted in \cite{UN} needed an extra term to
adjust the convergence when $q\uparrow1$ and gave no nice values at
negative integers. The choices $t=s-2,s-3,s-4,\ldots$ seem as good as
the value $\z(-m)$ is concerned, but extra poles at $s=2,3,4,\ldots$
emerge.  However, these poles disappear at the limit $q\uparrow1$. 
For example, with $t=s-2$ the residue at the simple pole $s=2$ is
$-(1-q)^2/\log q$  which goes to $0$ as $q\uparrow1$.  How things
become different depending on the choice of $t$ still seems to be
mysterious.

\noindent 2) If we introduce the $q$-analogue $\wt\z_q(s)$ of 
the alternating $\wt\z(s)$ in the introduction by
$$
\wt\z_q(s)=\sum_{n=1}^\infty (-1)^{n-1}\frac{q^{n(s-1)}}{[n]_q^{s}},
$$
the identity corresponding to \eqref{eq:alt} takes the form
$$
\wt\z_q(s) =\zeta_q(s)-2(1+q)^{-s}\zeta_{q^2}(s).
$$
In contrast to the situation of Euler, this does not help much and
indeed even makes things worse  because of the occurrence of another
base $q^2$.  It may be said that once $q$ is introduced, the
acceleration of convergence is fully achieved and nothing more is
needed.

The formula in Proposition 1 when specialized to $t=s-1$ becomes
\begin{eqnarray} \label{eq:zqexp}
\z_q(s)&=&(1-q)^{s}\sum_{r=0}^\infty \binom{s+r-1}{r}
\frac{q^{s+r-1}}{1-q^{s+r-1}}\\
&=& (1-q)^s\left(\frac{q^{s-1}}{1-q^{s-1}}+s\frac{q^{s}}{1-q^{s}}
+\frac{s(s+1)}2\frac{q^{s+1}}{1-q^{s+1}}
+\cdots\right). \nonumber
\end{eqnarray}

\noindent {\bf Proposition 2.} \begin{em}{\rm 1)} The function $\z_q(s)$ has 
  a simple pole at points in $1+2\pi i\Z/\!\log q$ and in the set
  $\left\{a+2\pi i b/\!\log q\, \vert a,b\in\Z,\,a\le0,\,b\ne0\right\}$. In
  particular, $s=1$ is a simple pole of $\z_q(s)$ with residue
  $(q-1)/\!\log q$.

\noindent {\rm 2)} For $m\in\Z,\ m\ge0$, the limiting 
value $\lim_{s\ra -m}\z_q(s)$ exists $($which we write $\z_q(-m))$ and
is given explicitly by
\begin{equation}
  \label{eq:zqneg}
  \z_q(-m)=(1-q)^{-m}\left\{
\sum_{r=0}^m(-1)^{r}\binom{m}{r}\frac1{q^{m+1-r}-1}
+\frac{(-1)^{m+1}}{(m+1)\log q}\right\}. 
\end{equation}
\end{em}

\noindent {\it Proof.}~Assertion 1) is straightforward from \eqref{eq:zqexp}, 
 the formula $\lim_{y\ra0}y/(1-q^y)=-1/\log q$ being used for the 
residue at $s=1$.  For 2), note the terms
with $r\ge m+2$ in the sum vanishes when $s\ra -m$ since $\binom{-m+r-1}{r}=0$
and $1-q^{-m+r-1}\ne0$.
On the other hand, for $r=m+1$ we have $\lim_{s\ra -m}(s+m)/(1-q^{s+m})
=-1/\log q$ and hence 
$$ \lim_{s\ra -m}\binom{s+m}{m+1}\frac{q^{s+m}}{1-q^{s+m}}
=\frac{(-1)^mm!}{(m+1)!}\left(-\frac1{\log q}\right)=\frac{(-1)^{m+1}}
{(m+1)\log q}.$$ The rest of the computation is clear. \qed

Before giving our general formula for $\lim_{q\uparrow1}\z_q(-m)$
(with expected value), let us look at the first few examples.

\noindent
{\it Example 1.} As stated in Proposition 2, $\zeta_q(s)$ has  a
simple pole at $s=1$ with residue $(q-1)/\log q$, which converges to
$1$ as $q\ra1$. This agrees with the well-known fact (reviewed later)
that $\zeta(s)$ has a simple pole at  $s=1$ with residue 1.
 
\noindent
{\it Example 2.} By \eqref{eq:zqneg} we have
$$ \z_q(0)=\frac1{q-1}-\frac1{\log q}.$$ Since
$$
\frac1{\log
  q}=\frac1{\log(1+(q-1))}=\frac1{(q-1)-(q-1)^2/2+\cdots}
=\frac1{q-1}+\frac12+O(q-1),
$$
we find $$ \lim_{q\ra1}\z_q(0)=-\frac12.$$
This agrees with Euler's computation $\zeta(0)=-1/2$.

\noindent
{\it Example 3.} Again by \eqref{eq:zqneg} we have
\begin{eqnarray*} \z_q(-1)&=&(1-q)^{-1}\left(\frac1{q^2-1}-\frac1{q-1}
+\frac1{2\log q}\right)\\
&=& \frac1{1-q}\left(\frac1{q-1}\cdot\frac1{2+q-1}-\frac1{q-1}
+\frac1{2\log q}\right)\\
&=&\frac1{1-q}\left(\frac1{2(q-1)}-\frac14+\frac{q-1}8+\cdots-\frac1{q-1}
+\frac1{2(q-1)}+\frac14-\frac{q-1}{24}+\cdots\right)\\
&&\!\!\!\!\!\!\!\!\!\!\longrightarrow -\frac1{12}\quad \text{as }q\ra1,
\end{eqnarray*}
in accordance with $\z(-1)=-1/12$.

Let the Bernoulli numbers $B_k$ be defined by the generating
series
$$
\frac{te^t}{e^t-1}\,(=\frac{t}{1-e^{-t}}) 
=\sum_{k=0}^\infty B_k\frac{t^k}{k!}.$$
First several values are
$$B_0=1,\ B_1=\frac12,\ B_2=\frac16,\ B_3=0,\ B_4=-\frac1{30},\
B_5=0,\ B_6=\frac1{42},\ B_7=0,\,\ldots$$
Now the
general formula is the following

\noindent {\bf Theorem 1.} \begin{em}For each non-negative integer $m$, 
  we have
  $$
  \lim_{q\uparrow1}\z_q(-m)=-\frac{B_{m+1}}{m+1}.
$$
\end{em}
\noindent {\it Proof.} ~On account of formula \eqref{eq:zqneg}, 
we have to show
$$\lim_{q\ra1}(1-q)^{-m}\left\{
  \sum_{r=0}^m(-1)^{r}\binom{m}{r}\frac1{q^{m+1-r}-1}
  +\frac{(-1)^{m+1}}{(m+1)\log q}\right\}=-\frac{B_{m+1}}{m+1}.$$
(Note here that since the sum on the left is finite so we may replace
the limit $q\uparrow1$ by $q\ra1$.) Multiplying both sides by
$(-1)^{m+1}(m+1)$ and changing $r\ra m+1-r$,  we see this is
equivalent to
$$\lim_{q\ra1}(1-q)^{-m}\left\{(m+1)
  \sum_{r=1}^{m+1}(-1)^{r}\binom{m}{r-1}\frac1{q^{r}-1}+\frac1{\log
    q}\right\}=(-1)^{m}B_{m+1}.$$
Writing 
$$\frac1{q^r-1}=\frac1r\cdot\frac{r\log q}{e^{r\log
q}-1}\cdot \frac1{\log q}$$ and using
$$\frac{t}{e^t-1}=\sum_{k=0}^\infty(-1)^k B_k\frac{t^k}{k!},$$ we have
\begin{eqnarray*}
&&(m+1)\sum_{r=1}^{m+1}(-1)^{r}\binom{m}{r-1}\frac1{q^r-1}\\
&=&(m+1)\sum_{r=1}^{m+1}(-1)^{r}\binom{m}{r-1}\frac1r\sum_{k=0}^\infty
(-1)^k B_k\frac{(r\log q)^{k}}{k!}\frac1{\log q}\\
&=&\sum_{k=0}^\infty\left(\sum_{r=1}^{m+1}(-1)^{r}\binom{m+1}{r}r^k\right)
(-1)^kB_k\frac{(\log q)^{k-1}}{k!}.
\end{eqnarray*}
Since the inner sum on the right can be calculated as 
\begin{eqnarray*}
\sum_{r=1}^{m+1}(-1)^{r}\binom{m+1}{r}r^k
&=&\left(\left(x\frac{d}{dx}\right)^k((1-x)^{m+1}-1)\right)\Bigg|_{x=1}\\
&=&\cases -1\; &\text{if}\quad k=0,\\
0\;&\text{if}\quad 0<k<m+1,\\
(-1)^{m+1}(m+1)!\;&\text{if}\quad k=m+1,
\endcases
\end{eqnarray*}
we find 
\begin{eqnarray*}
(m+1)\sum_{r=1}^{m+1}(-1)^{r}\binom{m}{r-1}\frac1{q^r-1}
=
-\frac1{\log q}+B_{m+1}(\log q)^{m}+ O((\log q)^{m+1}) \quad (\text{as} \,\;
q\to 1).
\end{eqnarray*}
 From this and the expansion  $\log q= q-1 + O((q-1)^2)\; (q\to1)$, we
obtain the desired result.\qed

\noindent {\it Remark.} In view of Theorem 1, we may define the
$q$-{\it Bernoulli number} $B_m(q)$ by 
$$B_m(q):=-m\z_q(1-m)\quad (m\ge1).$$
By \eqref{eq:zqneg} (letting
$m\to m-1$  and $r\to m-r$) we have the closed formula
\begin{eqnarray*}
  B_m(q)&=&(q-1)^{-m+1}\left\{\sum_{r=1}^m(-1)^r\binom{m}{r}\frac{r}{q^r-1}
+\frac1{\log q}\right\}\\
&=&(q-1)^{-m+1}\sum_{r=0}^m(-1)^r\binom{m}{r}\frac{r}{q^r-1}\quad (m\ge1).
\end{eqnarray*}
Here, the term with $r=0$ in the last sum should be read as $1/\log q$
(the limiting value when $r\to0$). This suggests to put
$$B_0(q)=\frac{q-1}{\log q}.$$ With this, the $q$-Bernoulli numbers
$\left\{B_m(q)\right\}_{m\ge0}$ satisfy the recursion
$$\sum_{m=0}^n(-1)^m\binom{n}{m}q^mB_m(q)=(-1)^nB_n(q)+\delta_{1n}\quad
(n\ge0),$$ where $\delta_{1n}=1$ if $n=1$ and $0$ otherwise, and the
generating function
$$F_q(t):=\sum_{m=0}^\infty B_m(q)\frac{t^m}{m!}$$ satisfies
the relation
$$F_q(qt)=e^tF_q(t)-te^t.$$ This $q$-Bernoulli number is essentially
(i.e., up to the sign $(-1)^m$) the same as the one introduced in
Tsumura \cite{T}.
 
The following fundamental relation, apart from its own importance,
guarantees that our computation at negative integers above does  give
us the correct values which we intended to obtain on a rigorous basis.

\noindent {\bf Theorem 2.} \begin{em}For any $s\in\C$, $s\ne1$, we have
$$
\lim_{q\uparrow1} \zeta_q(s)=\zeta(s).
$$\end{em}What we understand as the right-hand side  for arbitrary $s$
is the value of the function analytically continued to the whole
$s$-plane.  (We give the analytic continuation by using the {\it
  Euler-{\small Maclaurin} summation formula}, see the proof below.)
On the left-hand side, $q$ should  avoid the values with which
$\z_q(s)$ has a pole at $s$, but this is achieved  once $q$ gets close
enough to $1$.

\noindent {\it Example.} ~We give some numerical examples. Take $s=1/2$
and $q=0.999$ in \eqref{eq:zqexp}. Sum of the first $10^5$ terms gives
us the value $-1.46014527395\cdots$. Replacing $q$ by $q=0.99999$ and
taking the first $10^7$ terms we get the value $-1.460352417\cdots$,
which agrees with the actual value $\z(1/2)= -1.4603545088\cdots$ up
to $5$ decimal points. Take the point $s=1/2+14.1347i$ near the first
non-trivial zero ($=1/2+14.134725141734693790457251983562\cdots i$) of
$\z(s)$. For $q=0.9999$, the first $10^5$ terms gives the absurdly
large $10835.552\cdots  + 10270.785\cdots i$, while $10^6$ terms gives
$-0.000306477\cdots +  0.000794677\cdots i$  (the actual value is
$\z(1/2+14.1347i)=0.000003135364\cdots - 0.00001969336\cdots i$).  If
we take  $s=1/2+14.134725i$ and $q=0.99999$, the first $2\cdot 10^6 $
terms gives $-0.4690527\cdots - 0.4669811\cdots i$ and the $5\cdot
10^6$ terms $-0.000031064\cdots + 0.0000812513\cdots i$ (the actual
value is  $\z(1/2+14.134725i)= 0.000000017674\cdots -
0.00000011102\cdots i$).

Combining Theorem 1 and Theorem 2, we readily obtain

\noindent {\bf Corollary.} \begin{em}For each non-negative integer $m$, 
we have \end{em} $$
``1^m+2^m+3^m+4^m+5^m+\cdots"=\zeta(-m)=-\frac{B_{m+1}}{m+1}.
$$

\noindent
{\it Remarks.} 1)\ We can also define a $q$-analogue of the Hurwitz zeta
function $\z(s;a)=\sum_{n=0}^\infty 1/(n+a)^s$ by
$$\z_q(s;a)=\sum_{n=0}^\infty\frac{q^{(n+a)(s-1)}}{[n+a]_q^s}$$
and prove the identity
$$\lim_{q\uparrow1}\z_q(s;a)=\z(s;a)$$
for any $s\ne1$, as well as the formula 
$$\lim_{q\uparrow1}\z_q(-m;a)=-\frac{B_{m+1}(a)}{m+1}$$
for integers $m\le0$. Here, the Bernoulli polynomial $B_k(x)$ is defined by
the generating series
\begin{equation}
  \label{eq:bernpoly}
  \frac{te^{xt}}{e^t-1}=\sum_{k=0}^\infty B_k(x)\frac{t^k}{k!}.
\end{equation}
As in
the remark after Theorem 1, we can also define the $q$-Bernoulli
polynomial similarly and derive elementary formulas.  But to make our
presentation as concise as possible, we restrict ourselves to the case
of the Riemann zeta function.

\noindent 2)\  It would be amusing to note that the limit
\begin{equation}
  \label{eq:eisen}
  \lim_{q\uparrow1}(1-q)^k\sum_{n=1}^\infty
\frac{n^{k-1}q^n}{1-q^n}=(k-1)!\z(k)\quad (\forall k\ge2,\ k\in\Z)
\end{equation}
is derived easily from 
$$
\lim_{q\uparrow1}\z_q(k)=\zeta(k).$$
(The latter directly follows from
the definition  without appealing to Theorem 2 because we are in the
region of absolute convergence.) In fact, if we put $s=2$ in
\eqref{eq:zqexp} and make $r+1\ra n$, we have
$$\z_q(2)=(1-q)^2\sum_{n=1}^\infty \frac{nq^n}{1-q^n},$$
which gives the desired limit for $k=2$. For general $k$, we similarly
put $s=k$ in \eqref{eq:zqexp} and make $k+r-1\ra n$ to find
$$\z_q(k)=(1-q)^k\sum_{n=1}^\infty
\binom{n}{k-1}\frac{q^n}{1-q^n}.$$
(Observe $\binom{n}{k-1}=0$ for $n=1,2,\ldots,k-2$.) We note that
$$\binom{n}{k-1}=\frac{n^{k-1}}{(k-1)!}+\text{lower degree terms}$$
and, on taking the limit $q\uparrow1$, sums coming from lower terms
vanish inductively, hence we obtain the conclusion.

When $k$ is even and $k\ge4$, the series $\sum_{n=1}^\infty
n^{k-1}q^n/(1-q^n)=\sum_{n=1}^\infty
\left(\sum_{d|n}d^{k-1}\right)q^n$ constitutes the Fourier series of
the Eisenstein series $G_k(\tau)$ of weight $k$ on the modular group,
with constant term $-B_k/2k\,(=\z(1-k)/2)$. Here $\tau$ is a variable
in  the upper-half plane and is linked with $q$ by $q=e^{2\pi i
  \tau}$.  The modularity amounts to the transformation formula
$G_k(-1/\tau)=\tau^kG_k(\tau)$, which can be derived from, as Hecke
 \cite{H} showed, the functional equation of the corresponding
Dirichlet series $\varphi(s):=\z(s)\z(s+1-k)$:
$$(2\pi)^{-s}\G(s)\varphi(s)=(-1)^{k/2}(2\pi)^{s-k}\G(k-s)\varphi(k-s).$$
(When $k$ is odd, the functional equation of $\varphi(s)$ fails to
take this form and so the series $\sum_{n=1}^\infty
n^{k-1}q^n/(1-q^n)$ cannot be a Fourier series of a modular form.)
Hecke also showed that the residue of $\varphi(s)$ at the simple pole
$s=k$ is equal to  $(2\pi i)^kc_0/(k-1)!$ where $c_0$ is the constant
term of the corresponding modular form. In our case, the residue is
$\z(k)$ and thus the constant term of $G_k(\tau)$ is
$(k-1)!\z(k)/(2\pi i)^k=-B_k/2k$, as expected.  As an alternative way,
we may use \eqref{eq:eisen} to determine the constant term as follows:
Put $\tau=it$ with $t>0$. Then $e^{2\pi i(-1/it)}\ra0$ as $t\ra0$ and
so
\begin{eqnarray*}
  c_0&=&\lim_{t\ra 0}G_k(-\frac1{it})=\lim_{t\ra 0}(it)^kG_k(it)=
\lim_{q\uparrow1}\frac{(it)^k}{(1-q)^k}(1-q)^kG_k(it)\\
&=&\frac1{(2\pi i)^k}(k-1)!\zeta(k).
\end{eqnarray*}

\noindent{\it Proof of Theorem 2.} ~Recall the celebrated {\it summation 
formula of Euler} \cite{E1,E2} (and Maclaurin, cf.~\cite[\S7.21]{WW},
obtained simply by repeating integration by parts): For a
$C^\infty$-function $f(x)$ on $[1,\infty)$ and arbitrary integers
$M\ge0$, $N\ge1$, we have
\begin{eqnarray}
  \label{eq:sum}
  \sum_{n=1}^N f(n) &=& \int_1^Nf(x)\,dx+\frac12(f(1)+f(N))
+\sum_{k=1}^{M} \frac{B_{k+1}}{(k+1)!}\left(f^{(k)}(N)-f^{(k)}(1)\right)\\
&-&\frac{(-1)^{M+1}}{(M+1)!}\int_1^N\wt B_{M+1}(x)f^{(M+1)}(x)\,dx,\nonumber
\end{eqnarray}
where $\wt B_{M+1}(x)$ is the ``periodic Bernoulli polynomial'' defined by 
$$\wt B_k(x)=B_k(x-[x])\quad ([x] \text{ is the largest integer not
  exceeds } x).$$
Recall the Bernoulli polynomial $B_k(x)$ is defined by the generating
series \eqref{eq:bernpoly}:
$$B_0(x)=1,\ B_1(x)=x-\frac12,\ B_2(x)=x^2-x+\frac16,\ B_3(x)=x^3-\frac32
x^2+\frac12 x,\,\ldots$$
As is well-known, by taking $f(x)=x^{-s}$ and
letting $N\ra\infty$, we obtain the analytic continuation of $\z(s)$
to the region $\re(s)>-M$:
\begin{equation}
  \label{eq:riemann}
  \z(s)=\frac1{s-1}+\frac12+\sum_{k=1}^{M} \frac{B_{k+1}}{(k+1)!}(s)_k
-\frac{(s)_{M+1}}{(M+1)!}\int_1^\infty \wt B_{M+1}(x)x^{-s-M-1}\,dx,
\end{equation}
where $(s)_k:=s(s+1)\cdots(s+k-1)$. Since we may choose $M$ arbitrary
large, this gives the analytic  continuation of $\z(s)$ to the whole
$s$-plane, revealing the (unique) simple pole at $s=1$ with residue
$1$.

Now we take $f(x)=q^{x(s-1)}/(1-q^x)^s$ and $M=1$ in \eqref{eq:sum}. 
Assuming  $\re(s)>1$ and noting 
\begin{eqnarray*}
  f'(x)&=&\log q\cdot q^{x(s-1)}\frac{s-1+q^x}{(1-q^x)^{s+1}},\\
f''(x)&=&(\log q)^2 q^{x(s-1)}\frac{s(s+1)-
3s(1-q^x)+(1-q^x)^2}{(1-q^x)^{s+2}},
\end{eqnarray*} and in general $f^{(k)}(x)=(\log q)^kq^{x(s-1)}(1-q^x)^{-s-k}
\times(\text{a polynomial in }s\text{ and }q^x)$, we see that we can
take the limit $N\ra\infty$ and obtain
\begin{eqnarray*}
  \sum_{n=1}^\infty\frac{q^{n(s-1)}}{(1-q^n)^s}&=&
\int_1^\infty\frac{q^{x(s-1)}}{(1-q^x)^s}\,dx+\frac12\cdot
\frac{q^{s-1}}{(1-q)^s}
-\frac1{12}(\log q)q^{s-1}\frac{s-1+q}{(1-q)^{s+1}}\\
&-&\frac{(\log q)^2}2\int_1^\infty \wt B_2(x)q^{x(s-1)}\frac{s(s+1)-
3s(1-q^x)+(1-q^x)^2}{(1-q^x)^{s+2}}\,dx
\end{eqnarray*} for $\re(s)>1$. 
The first integral on the right is computed as
$$\int_1^\infty\frac{q^{x(s-1)}}{(1-q^x)^s}\,dx=
\int_1^\infty\frac{q^{-x}}{(q^{-x}-1)^s}\,dx
=\left[\frac{(q^{-x}-1)^{1-s}}{(s-1)\log q}\right]_1^\infty
=-\frac{q^{s-1}(1-q)^{1-s}}{(s-1)\log q}.$$
We therefore obtain
\begin{eqnarray}
  \label{eq:zqeul}
 \z_q(s) &=&(1-q)^s\sum_{n=1}^\infty\frac{q^{n(s-1)}}{(1-q^n)^s}\\
&=&\frac{q^{s-1}}{s-1}\cdot \frac{q-1}{\log q}+\frac{q^{s-1}}2
+\frac{q^{s-1}}{12}\frac{\log q}{q-1}(s-1+q)\nonumber \\
&-&(1-q)^s\frac{(\log q)^2}2\int_1^\infty\!\wt B_2(x)q^{x(s-1)}
\frac{s(s+1)-
3s(1-q^x)+(1-q^x)^2}{(1-q^x)^{s+2}}\,dx.\nonumber 
\end{eqnarray}
Unlike the classical case \eqref{eq:riemann}, just to let $M$ larger
does not make the convergence of the integral better, since the factor
$q^{x(s-1)}$ in $f^{(M+1)}(x)$ always forces $\re(s)>1$.  Instead, we
use in \eqref{eq:zqeul} the Fourier expansion of the periodic
Bernoulli polynomials\footnote{We owe Ueno-Nishizawa \cite{UN} the
  idea of replacing $\wt B_2(x)$ in the integral by its  Fourier
  expansion. However, our argument that follows, which uses only
  integration by parts and no confluent hypergeometric functions or
  the like, seems considerably different  from the one in \cite{UN}.}
(cf.~\cite[Ch.IX, Misc.~Ex.~12]{WW})
\begin{equation}
  \label{eq:fourier}
  \wt B_k(x)=-k!\sum_{n\in\Z\atop n\ne0}
\frac{e^{2\pi i nx}}{(2\pi in)^k}.
\end{equation}
The equality is valid for all real numbers $x$ when $k\ge2$, the sum
being absolutely and uniformly convergent.  Putting this (for $k=2$)
into \eqref{eq:zqeul} and interchanging the summation and the
integration, we find
\begin{eqnarray*}
  \z_q(s) &=&\frac{q^{s-1}}{s-1}\cdot \frac{q-1}{\log q}+\frac{q^{s-1}}2
+\frac{q^{s-1}}{12}\frac{\log q}{q-1}(s-1+q)\\
&+&(1-q)^s(\log q)^2\sum_{n\in\Z\atop n\ne0}
\frac1{(2\pi in)^2}\int_1^\infty e^{2\pi i nx}q^{x(s-1)}
\frac{s(s+1)-3s(1-q^x)+(1-q^x)^2}{(1-q^x)^{s+2}}\,dx.
\end{eqnarray*}
Further we make a change of variable $q^x=u$ to obtain
\begin{eqnarray} \label{eq:zqbq}
  \z_q(s) &=&\frac{q^{s-1}}{s-1}\cdot \frac{q-1}{\log q}+\frac{q^{s-1}}2
+\frac{q^{s-1}}{12}\frac{\log q}{q-1}(s-1+q)\\ \nonumber
&&-(1-q)^s\log q\sum_{n\in\Z\atop n\ne0}
\frac1{(2\pi in)^2}\left\{s(s+1)b_q(s-1+\de n,-s-1)\right.\\ \nonumber
&&\left.\qquad -3s b_q(s-1+\de n,-s)+b_q(s-1+\de n,-s+1)\right\},
\end{eqnarray}
where we have put $\de=2\pi i/\log q$  and 
$$b_t(\alpha,\beta)=\int_0^t u^{\alpha-1}(1-u)^{\beta-1}\,du$$
(referred to as the incomplete beta function\footnote{We are tempted
  to  remind the reader that the beta integral is often called the
  {\it Euler integral}.}).  Note that each of the  incomplete beta
integrals in \eqref{eq:zqbq} converges absolutely for $\re(s)>1$ and
uniformly bounded with respect to $n$;
$$\left\vert b_q(s-1+\de n,-s+\nu)\right\vert\le\int_0^qu^{\sigma-2}
(1-u)^{-\sigma+\nu-1}du\quad(\forall
n,\,\sigma=\re(s),\,\nu=-1,0,1),$$
hence the sum converges absolutely.

Now, repeated use of integration by parts provides us the formula
\begin{eqnarray*}
b_t(\alpha,\beta)&=&\int_0^t\left(\frac{u^\alpha}{\alpha}\right)'
(1-u)^{\beta-1}\,du=\frac1{\alpha}t^\alpha(1-t)^{\beta-1}-
\frac{1-\beta}{\alpha}\int_0^tu^\alpha(1-u)^{\beta-2}\,du\\
&=&\frac1{\alpha}t^\alpha(1-t)^{\beta-1}-
\frac{1-\beta}{\alpha}\int_0^t\left(\frac{u^{\alpha+1}}{\alpha+1}\right)'
(1-u)^{\beta-2}\,du\\ \noalign{\vskip1ex}
&=&\cdots\cdots\\ \noalign{\vskip1ex}
&=&\sum_{k=1}^{M-1}(-1)^{k-1}\frac{(1-\beta)_{k-1}}{(\alpha)_k}
t^{\alpha+k-1}(1-t)^{\beta-k}\\
&&\qquad \qquad \qquad
+(-1)^{M-1}\frac{(1-\beta)_{M-1}}{(\alpha)_{M-1}}\beta_t(\alpha+M-1,\beta-M+1)
\end{eqnarray*}
for any
$M\ge2$. Applying this to $b_q(s-1+\de n,-s-1)$, we have (note $q^{\de
  n}=1$)
\begin{eqnarray*}
 b_q(s-1+\de n,-s-1)&=&\sum_{k=1}^{M-1}(-1)^{k-1}\frac{(s+2)_{k-1}}
{(s-1+\de n)_k}
q^{s+k-2}(1-q)^{-s-1-k}\\
&&+(-1)^{M-1}\frac{(s+2)_{M-1}}{(s-1+\de n)_{M-1}}
b_q(s-2+M+\de n,-s-M).
\end{eqnarray*}
This accomplishes the analytic continuation of $b_q(s-1+\de n,-s-1)$
as a function of $s$ into the region $\re(s)>2-M$.  From this we have
\begin{eqnarray*}
&&  \sum_{n\in\Z\atop n\ne0}\frac{s(s+1)}{(2\pi in)^2}b_q(s-1+\de n,-s-1)\\
&=&\sum_{k=1}^{M-1}(-1)^{k-1}\sum_{n\in\Z\atop n\ne0}
\frac1{(2\pi in)^2}\frac{(s)_{k+1}}{(s-1+\de n)_k}q^{s+k-2}(1-q)^{-s-1-k}\\
&+&(-1)^{M-1}\sum_{n\in\Z\atop n\ne0}\frac1{(2\pi in)^2}
\frac{(s)_{M+1}}{(s-1+\de n)_{M-1}}\int_0^q u^{s-3+M+\de n}(1-u)^{-s-M-1}\,du\\
&=&\sum_{k=1}^{M-1}(-1)^{k-1}\sum_{n\in\Z\atop n\ne0}
\frac1{(2\pi in)^2}\frac{(s)_{k+1}}{(s-1+\de n)_k}q^{s+k-2}(1-q)^{-s-1-k}\\
&-&(-1)^{M-1}\log q\sum_{n\in\Z\atop n\ne0}\frac1{(2\pi in)^2}
\frac{(s)_{M+1}}{(s-1+\de n)_{M-1}}\int_1^\infty e^{2\pi i nx}
q^{x(s-2+M)}(1-q^x)^{-s-M-1}\,dx.
\end{eqnarray*}
Using $$\lim_{q\to1}\frac{\log q}{1-q}=-1,\quad
\lim_{q\to1}(1-q)^k(s-1+\de n)_k=(-2\pi i)^k,\quad
\lim_{q\to1}\frac{1-q^x}{1-q}=x,$$ we obtain, for $\re(s)>2-M$,
\begin{eqnarray*}
&&  \lim_{q\uparrow1}
(1-q)^s\log q\sum_{n\in\Z\atop n\ne0}\frac{s(s+1)}{(2\pi in)^2}
b_q(s-1+\de n,-s-1)\\ &=&
\sum_{k=1}^{M-1}\sum_{n\in\Z\atop n\ne0}\frac{(s)_{k+1}}{(2\pi i n)^{k+2}}
-\sum_{n\in\Z\atop n\ne0}\frac{(s)_{M+1}}{(2\pi i n)^{M+1}}
\int_1^\infty e^{2\pi i nx}x^{-s-M-1}\,dx\\
&=& -\sum_{k=1}^{M-1}\frac{B_{k+2}}{(k+2)!}(s)_{k+1}+\frac{(s)_{M+1}}{(M+1)!}
\int_1^\infty \wt B_{M+1}(x)x^{-s-M-1}\,dx.
\end{eqnarray*}
In the last equality, we have used \eqref{eq:fourier} and its
specialization ($x=1$)
$$\sum_{n\in\Z\atop n\ne0}\frac1{(2\pi i n)^k}=-\frac{B_k}{k!}$$
valid
for all $k\ge2$.  We do exactly the same for the terms containing
$b_q(s-1+\de n,-s)$ and  $b_q(s-1+\de n,-s+1)$. As it turns out
however, the contributions from these two  vanish when we take
$q\uparrow1$, for the powers of $1-q$ involved  are lower  than those
from $b_q(s-1+\de n,-s-1)$. We therefore obtain, for $\re(s)>2-M$,
\begin{eqnarray*}
  \lim_{q\uparrow1}\z_q(s)&=&\frac1{s-1}+\frac12+\frac{s}{12}
+\sum_{k=2}^{M}\frac{B_{k+1}}{(k+1)!}(s)_{k}-\frac{(s)_{M+1}}{(M+1)!}
\int_1^\infty \wt B_{M+1}(x)x^{-s-M-1}\,dx\\
&=&\frac1{s-1}+\frac12+\sum_{k=1}^{M}\frac{B_{k+1}}{(k+1)!}(s)_{k}
-\frac{(s)_{M+1}}{(M+1)!}
\int_1^\infty \wt B_{M+1}(x)x^{-s-M-1}\,dx.
\end{eqnarray*}
This coincides with formula \eqref{eq:riemann} for $\z(s)$
valid in $\re(s)>-M$,  and thus the theorem is established since the
integer $M$ can be arbitrary large.  \qed



\vspace{3ex}
\noindent Faculty of Mathematics, Kyushu University, \\
Hakozaki, Fukuoka, 812-8581 JAPAN\\
{\tt mkaneko@math.kyushu-u.ac.jp}\\
\\
Department of Mathematics, 
Tokyo Institute of Technology, \\
Meguro, Tokyo, 152-0033 JAPAN\\
{\tt kurokawa@math.titech.ac.jp}\\
\\
Faculty of Mathematics, Kyushu University, \\
Hakozaki, Fukuoka, 812-8581 JAPAN\\
{\tt wakayama@math.kyushu-u.ac.jp}

\end{document}